\documentclass[twoside,emlines,bezier,12pt]{article}
\usepackage{amsfonts,amssymb}

\begin{document}

\begin{center}

{\Large \bf $k$-Congruences on semirings}

$ $

\bigskip {Song-Chol Han}

\vskip 3mm

{\small \centerline{{\it Faculty of Mathematics, Kim Il Sung
University, Pyongyang, Democratic People's Republic of Korea}}}

\renewcommand{\thefootnote}{\alph{footnote}}
\setcounter{footnote}{-1} \footnote{{\it E-mail address:}
ryongnam5@yahoo.com.}
\end{center}

\vspace{0.3cm}

{\small \noindent {\bf ABSTRACT}
\medskip

\noindent For any semiring, the concept of $k$-congruences is
introduced, criteria for $k$-congruences are established, it is
proved that there is an inclusion-preserving bijection between
$k$-congruences and $k$-ideals, and an equivalent condition for
the existence of a zero is presented with the help of
$k$-congruences. It is shown that a semiring is $k$-simple iff it
is $k$-congruence-simple, and that inclines are $k$-simple iff
they have at most 2 elements. Lemma 2.12(i) in [Glas. Mat. 42(62)
(2007) 301] is pointed out being false.
\medskip

\medskip

}
\medskip

$ $

{\section*{\noindent \bf 1. Introduction}}

$ $

The notion of semirings was introduced by Vandiver
\cite{Vandiver1934} in 1934.

The most trivial example of a semiring which is not a ring is the
set of all nonnegative integers with the usual addition and
multiplication that is the first algebraic structure we encounter
in life. The semiring of all ideals of a commutative ring with
identity and the semiring of all endomorphisms on a commutative
semigroup are nontrivial examples of semirings.

Semirings are algebraic systems that generalize both rings and
distributive lattices and have a wide range of applications
recently in diverse branches of mathematics and computer science
such as graph, optimization, automata, formal language, algorithm,
coding theory, cryptography, etc
\cite{Golan1999,Hebisch1998,Maze2007,Sakalauskas2005}.

Semirings have two binary operations of addition and
multiplication which are connected by the ring-like distributive
laws. But unlike in rings, subtraction is not allowed in
semirings. For this reason, there are considerable differences
between ring theory and semiring theory and many results in ring
theory have no analogues in semirings. For instance, one can note
that for any ring, there is a one-to-one correspondence between
the congruences and ideals which associates with each congruence
its zero-class, thus a ring is congruence-simple iff it is
ideal-simple. However, the same statement is false for semirings.
While it is true that every ideal in a semiring induces a
congruence, there are semirings which have congruences not induced
by any ideal.

In order to narrow the gap, Henriksen \cite{Henriksen1958} in 1958
defined $k$-ideals in semirings, which are a special kind of
semiring ideals much closer to ring ideals than the general ones.
Since then, many researchers have developed $k$-ideal theory
\cite{Atani2007,Bhuniya2012,Han2015,LaTorre1965,Lescot2015,Sen1993,Weinert1996,Yesilot2010}.

The objective of the present paper is to study the connection of
$k$-ideals with congruences in semirings.

With this theme, Lescot \cite{Lescot2011} showed the existence of
a bijection between the saturated ideals and excellent congruences
in any characteristic one semiring, and Zhou and Yao
\cite{Zhou2011} proved the existence of a one-to-one
correspondence between the ideals which are lower sets and regular
congruences in any additively idempotent semiring. It turns out
that characteristic one semirings are additively idempotent
commutative semirings and saturated ideals are just $k$-ideals.
And it was proved in Han \cite{Han2015} that the ideals which are
lower sets are nothing but $k$-ideals in additively idempotent
semirings.

In this paper, for any semiring, the concept of $k$-congruences is
introduced, criteria for $k$-congruences are established (Theorem
3.3), it is proved that there is an inclusion-preserving bijection
between $k$-congruences and $k$-ideals (Theorem 3.8), and an
equivalent condition for the existence of a zero is presented with
the help of $k$-congruences (Theorem 4.2). It is shown that a
semiring is $k$-simple iff it is $k$-congruence-simple (Theorem
4.4), and that inclines are $k$-simple iff they have at most 2
elements (Theorem 5.4). Lemma 2.12(i) in \cite{Atani2007} is
pointed out to be false (Example 6.1).
\medskip

$ $

{\section*{\noindent \bf 2. Preliminaries}}

$ $

In this section, we recall some known definitions and facts
\cite{ElBashir2007,Golan1999,Han2015,Hebisch1998}.

Throughout this paper, $R$ denotes a {\it semiring} $(R,+,\cdot)$,
unless otherwise stated. That is to say, $(R,+)$ is a commutative
semigroup, $(R,\cdot)$ is a semigroup, and multiplication
distributes over addition from either side.

If $R$ has an additively neutral element $0$ and $0r=r0=0$ for all
$r\in R$, then $0$ is called a {\it zero} of $R$. We write $0\in
R$ when $R$ has a zero $0$. If $R$ has a multiplicatively neutral
element $1$, then $1$ is called an {\it identity} of $R$.

$R$ is called a {\it commutative} semiring if multiplication is
commutative.

A nonempty subset $A$ of $R$ is called an {\it ideal} of $R$ if
$a+b\in A$ and $ra,ar\in A$ for all $a,b\in A$ and $r\in R$. $R$
and $\{0\}$ (if $0\in R$) are said to be {\it trivial} ideals of
$R$. Denote by $\mathcal{I}(R)$ the family of all ideals of $R$.

For an ideal $A$ of $R$, the set
$$
\overline{A}=\{x\in R\mid x+a=b\ \mbox{for some}\ a,b\in A\}
$$
is called the {\it subtractive closure} or {\it $k$-closure} of
$A$ in $R$. Then $\overline{A}$ is an ideal of $R$ and it holds
that $A\subseteq \overline{A}$ and
$\overline{(\overline{A})}=\overline{A}$. $A$ is said to be {\it
$k$-closed} in $R$ if $A=\overline{A}$. For ideals $A$ and $B$ of
$R$, $A\subseteq B$ implies $\overline{A}\subseteq \overline{B}$.

An ideal $A$ of $R$ is called a {\it subtractive ideal} or {\it
$k$-ideal} of $R$ if $A=\overline{A}$. If $A$ is an ideal of $R$,
then $\overline{A}$ is a $k$-ideal of $R$. $R$ is a $k$-ideal of
itself, and $\{0\}$ is also a $k$-ideal of $R$ if $0\in R$. Denote
by $\mathcal{KI}(R)$ the family of all $k$-ideals of $R$.

An equivalence relation $\equiv$ on $R$ is called a {\it
congruence} on $R$ if for any $a,b,c\in R$, $a\equiv b$ implies
$a+c\equiv b+c$, $ac\equiv bc$ and $ca\equiv cb$. Denote by
$\mathcal{C}(R)$ the family of all congruences on $R$.

Given a congruence $\theta$ on $R$, the quotient set
$R/\theta=\{[x]\mid x\in R\}$ consisting of all the congruence
classes forms a semiring under the operations defined as
$[x]+[y]=[x+y]$ and $[x]\cdot [y]=[xy]$ for $x,y\in R$. This
semiring $R/\theta$ is called the {\it quotient semiring} of $R$
by $\theta$. If $0\in R$, then $[0]$ is a zero of $R/\theta$.

An ideal $A$ of $R$ defines a congruence $\kappa_A$ on $R$ by
\[
x \kappa_A y\ \Longleftrightarrow\ x+a=y+b\ \mbox{for some}\
a,b\in A.
\]
Then the $k$-closure $\overline{A}$ of $A$ is a zero of the
quotient semiring $R/\kappa_A$. In addition,
$\kappa_A=\kappa_{\overline{A}}$.

Below, we assume that $|R|\geqslant 2$ to avoid trivial
exceptions.
\medskip

$ $

{\section*{\noindent \bf 3. $k$-Congruences on semirings}}

$ $

In this section, we prove that there is a bijection between the
family of all $k$-congruences and the family of all $k$-ideals for
any semiring.
\medskip

\noindent {\bf Lemma 3.1.}\quad  If $\theta$ is a congruence on
$R$ and $R/\theta$ has a zero $0_\theta$, then $0_\theta$ is a
$k$-ideal of $R$.
\medskip

\noindent {\bf Proof.}\quad  Show that $0_\theta$ is an ideal of
$R$. If $a,b\in 0_\theta$, then
$[a+b]=[a]+[b]=0_\theta+0_\theta=0_\theta$ and so $a+b\in
0_\theta$. If $a\in 0_\theta$ and $r\in R$, then $[ar]=[a]\cdot
[r]=0_\theta\cdot [r]=0_\theta$ and similarly $[ra]=0_\theta$,
thus $ar,ra\in 0_\theta$. Show that $0_\theta$ is $k$-closed in
$R$. If $x+b\in 0_\theta$ and $b\in 0_\theta$, then
$[x]=[x]+0_\theta=[x]+[b]=[x+b]=0_\theta$ and so $x\in 0_\theta$.
\quad  $\square$
\medskip

\noindent {\bf Definition 3.2.}\quad  A congruence $\theta$ on $R$
is called a {\it subtractive congruence} or {\it $k$-congruence}
if there is an ideal $A$ of $R$ such that $\theta=\kappa_A$.
Denote by $\mathcal{KC}(R)$ the family of all $k$-congruences on
$R$.
\medskip

\noindent {\bf Theorem 3.3.}\quad  If $\theta$ is a congruence on
$R$, then the following conditions are equivalent.
\medskip

\noindent (1)\quad  $\theta$ is a $k$-congruence.

\noindent (2)\quad  $R/\theta$ has a zero $0_\theta$ and
$\theta\subseteq \kappa_{0_\theta}$.

\noindent (3)\quad  $R/\theta$ has a zero $0_\theta$ and
$\theta=\kappa_{0_\theta}$.
\medskip

\noindent {\bf Proof.}\quad  (1) $\Rightarrow$ (2)\quad Suppose
that $\theta=\kappa_A$ for an ideal $A$ of $R$. Then the
$k$-closure $\overline{A}$ is a zero of $R/\theta$. If $x\theta
y$, then $x+a=y+b$ for some $a,b\in A$ while $A\subseteq
\overline{A}$ and so $x\kappa_{\overline{A}} y$.

(2) $\Rightarrow$ (3)\quad  By Lemma 3.1, $0_\theta$ is an ideal
of $R$. Show that $\kappa_{0_\theta}\subseteq\theta$. If
$x\kappa_{0_\theta} y$, then $x+a=y+b$ for some $a,b\in 0_\theta$
and
$[x]=[x]+0_\theta=[x]+[a]=[x+a]=[y+b]=[y]+[b]=[y]+0_\theta=[y]$,
thus $x\theta y$. Hence $\theta=\kappa_{0_\theta}$.

(3) $\Rightarrow$ (1)\quad  This follows from Lemma 3.1.\quad
$\square$
\medskip

\noindent {\bf Example 3.4.}\quad  Let $R=\{0,a,b,c,d,1\}$ be the
lattice with the Hasse graph shown in Figure 1, which is not a
distributive lattice. Define a multiplication $\cdot$ on $R$ by
\[
x\cdot y=\left\{
\begin{array}{ccccc}
d, &\mbox{if}\  x,y\in \{1,b,c,d\}\\
0, &\mbox{otherwise}.
\end{array}
\right.
\]
Then $(R,\vee,\cdot)$ is a semiring and $A=\{0,a\}$ is an ideal of
$(R,\vee,\cdot)$. Thus $\kappa_A$ is a $k$-congruence on
$(R,\vee,\cdot)$ with congruence classes $\{0,a\}$ and
$\{1,b,c,d\}$.
\medskip

\begin{picture}(0,95)
\put(200,54){\circle{3}} \put(200,78){\circle{3}}
\put(176,54){\circle{3}} \put(224,54){\circle{3}}
\put(200,30){\circle{3}} \put(212,42){\circle{3}}
\put(200,55.5){\line(0,1){21}} \put(201,53){\line(1,-1){10}}
\put(199,77){\line(-1,-1){22}} \put(201,77){\line(1,-1){22}}
\put(201,31){\line(1,1){10}} \put(199,31){\line(-1,1){22}}
\put(213,43){\line(1,1){10}} \put(197,81){$\small 1$}
\put(168,52){$\small a$} \put(193,52){$\small b$}
\put(227,52){$\small c$} \put(212,35){$\small d$}
\put(197,19){$\small 0$}

\put(180,0){Figure 1}
\end{picture}
\medskip

\noindent {\bf Example 3.5.}\quad  Let $R=\{0,a,b,c,1\}$ be the
distributive lattice with the Hasse graph shown in Figure 2. Then
the congruence $\theta$ on $(R,\vee,\wedge)$ with congruence
classes $\{0\}$ and $\{1,a,b,c\}$ is not a $k$-congruence. In
fact, the zero of $R/\theta$ is $0_\theta=\{0\}$ and $a\theta b$
while $a\vee 0\neq b\vee 0$.
\medskip

\begin{picture}(0,80)
\put(200,30){\circle{3}} \put(200,45){\circle{3}}
\put(200,65){\circle{3}} \put(180,55){\circle{3}}
\put(220,55){\circle{3}} \put(201.5,45.5){\line(2,1){17}}
\put(198.5,45.5){\line(-2,1){17}}
\put(201.5,64.5){\line(2,-1){17}}
\put(198.5,64.5){\line(-2,-1){17}} \put(200,31.5){\line(0,1){12}}
\put(197,19){$\small 0$} \put(172,53){$\small a$}
\put(223,53){$\small b$} \put(203,40){$\small c$}
\put(197,68){$\small 1$}

\put(180,0){Figure 2}
\end{picture}
\medskip

By Lemma 3.1 and Theorem 3.3, we can define a function $\iota:
\mathcal{KC}(R)\rightarrow \mathcal{KI}(R)$ by
$\iota(\theta)=0_\theta$, where $0_\theta$ is a zero of
$R/\theta$. Also we define a function
$\kappa:\mathcal{I}(R)\rightarrow \mathcal{C}(R)$ by
$\kappa(A)=\kappa_A$.

Then we obtain the following conclusion from Theorem 3.3.
\medskip

\noindent {\bf Lemma 3.6}\quad  If $\theta$ is a $k$-congruence on
$R$, then $\kappa(\iota(\theta))=\theta$.
\medskip

\noindent {\bf Lemma 3.7.}\quad  An ideal $A$ of $R$ is a
$k$-ideal of $R$ iff $\iota(\kappa(A))=A.$
\medskip

\noindent {\bf Proof.}\quad  This follows from the fact that
$\iota(\kappa(A))=\iota(\kappa_A)=0_{\kappa_A}=\overline{A}$.\quad
$\square$
\medskip

\noindent {\bf Theorem 3.8.}\quad  (1)\quad  The restriction of
$\kappa$ to $\mathcal{KI}(R)$ is an inclusion-preserving bijection
from $\mathcal{KI}(R)$ onto $\mathcal{KC}(R)$.

\noindent (2)\quad  If $\mathcal{KI}(R)\subseteq
\mathcal{F}\subseteq \mathcal{I}(R)$ and $\kappa$ is injective on
$\mathcal{F}$, then $\mathcal{F}=\mathcal{KI}(R).$
\medskip

\noindent {\bf Proof.}\quad  (1)\quad  This follows from Lemmas
3.6 and 3.7.

(2)\quad  If $A\in \mathcal{F}$, then $\overline{A}\in
\mathcal{KI}(R)\subseteq \mathcal{F}$. Since
$\kappa(A)=\kappa_A=\kappa_{\overline{A}}=\kappa(\overline{A})$
and $\kappa$ is injective on $\mathcal{F}$, $A=\overline{A}$ and
thus $A\in \mathcal{KI}(R)$.\quad $\square$
\medskip

\noindent {\bf Remark 3.9.}\quad  Let $R$ be an additively
idempotent semiring, i.e. $r+r=r$ for all $r\in R$. If $A$ is an
ideal of $R$ and $x,y\in R$, then the condition that $x+a=y+b$ for
some $a,b\in A$ implies the condition that $x+c=y+c$ for some
$c\in A$. In fact, putting $c=a+b\in A$ gives
$x+c=x+a+b=x+a+a+b=y+b+a+b=y+a+b=y+c$.
\medskip

\noindent {\bf Remark 3.10.}\quad  A characteristic one semiring
is a commutative semiring with identity $1$ such that $1+1=1$
\cite{Lescot2011,Lescot2015}. Obviously, characteristic one
semirings are additively idempotent.
\medskip

Remarks 3.9 and 3.10 show that Theorem 3.8(1) is a generalization
of Theorem 3.8 in \cite{Lescot2011} and Theorem 5 in
\cite{Zhou2011}.
\medskip

$ $

{\section*{\noindent \bf 4. $k$-Simple semirings}}

$ $

In this section, we present an equivalent condition for a semiring
to have a zero and a necessary and sufficient condition for a
semiring to be $k$-simple by means of $k$-congruences.
\medskip

\noindent {\bf Remark 4.1}\quad  $R\times R$ is always a
$k$-congruence on $R$ because $\kappa_R=R\times R$. But
$\mbox{id}_R$ is not necessarily a $k$-congruence on $R$.
\medskip

\noindent {\bf Theorem 4.2.}\quad  $R$ has a zero $0$ iff
$\mbox{id}_R$ is a $k$-congruence on $R$. In this case,
$\kappa_{\{0\}}=\mbox{id}_R$.
\medskip

\noindent {\bf Proof.}\quad  (Necessity)\quad  If $R$ has a zero
$0$, then $\{0\}$ is an ideal of $R$ and obviously
$\kappa_{\{0\}}=\mbox{id}_R$.

(Sufficiency)\quad  If $\mbox{id}_R$ is a $k$-congruence on $R$,
then each of the congruence classes is a singleton. By Theorem
3.3, $R/{\mbox{id}_R}$ has a zero $[a]=\{a\}$, where $a\in R$.
Show that $a$ is a zero of $R$. For any $x\in R$,
$[x+a]=[x]+[a]=[x]$ and so $x+a=x$. Also $[xa]=[x]\cdot [a]=[a]$
and similarly $[ax]=[a]$, thus $xa=a=ax$.\quad $\square$
\medskip

Theorem 4.2 enables us to give the following definition.
\medskip

\noindent {\bf Definition 4.3.}\quad  $R$ is said to be {\it
$k$-congruence-simple} if it admits no $k$-congruences other than
$R\times R$ and $\mbox{id}_R$.
\medskip

$R$ is said to be {\it $k$-simple} or {\it simple} if it has no
nontrivial $k$-ideals \cite{Bhuniya2012,Stone1977}. $R$ is said to
be {\it congruence-simple} if it has just two congruences
\cite{ElBashir2007}. $R$ is said to be {\it ideal-free} if it has
no nontrivial ideals \cite{Stone1977}. Obviously,
congruence-simple semirings are $k$-congruence-simple and
ideal-free semirings are $k$-simple.
\medskip

\noindent {\bf Theorem 4.4.}\quad  $R$ is $k$-simple iff it is
$k$-congruence-simple.
\medskip

\noindent {\bf Proof.}\quad  By Remark 4.1 and Theorem 4.2, this
follows from the bijection established in Theorem 3.8(1).\quad
$\square$
\medskip

\noindent {\bf Corollary 4.5.}\quad  If $R$ is a ring, then the
following conditions are equivalent.
\medskip

\noindent (1) $R$ is ideal-simple as a ring.

\noindent (2) $R$ is ideal-free as a semiring.

\noindent (3) $R$ is $k$-simple as a semiring.

\noindent (4) $R$ is congruence-simple as a ring.

\noindent (5) $R$ is congruence-simple as a semiring.

\noindent (6) $R$ is $k$-congruence-simple as a semiring.
\medskip

\noindent {\bf Proof.}\quad  Proposition 1.2 in
\cite{ElBashir2007} shows that (1), (2), (4) and (5) are
equivalent. For any ring, a semiring ideal is a ring ideal iff it
is a $k$-ideal. Thus (1) and (3) are equivalent. Hence Theorem 4.4
completes the proof.\quad $\square$
\medskip

$ $

{\section*{\noindent \bf 5. $k$-Simple inclines}}

$ $

In this section, we show that inclines are $k$-simple iff they
have at most 2 elements.

On any idempotent commutative semigroup $(S,+)$, a partial
ordering $\leqslant$ is defined by $x\leqslant y\Leftrightarrow
x+y=y$ for $x,y\in S$. Then the poset $(S,\leqslant)$ is a
join-semilattice, where $x\vee y=x+y$ for all $x,y\in S$.

If $R$ is an additively idempotent semiring, then $x\leqslant y$
implies $zx\leqslant zy$ and $xz\leqslant yz$ for all $x,y,z\in
R$.

An additively idempotent semiring $R$ is called an {\it incline}
if $x+xy=x=x+yx$ for all $x,y\in R$. Then $xy\leqslant x$ and
$yx\leqslant x$ for all $x,y\in R$ \cite{Cao1984}.

The following are examples of inclines:

\noindent (1)\quad The Boolean algebra $(\{0,1\},\vee,\wedge)$,

\noindent (2)\quad Each distributive lattice $(D,\vee,\wedge)$,

\noindent (3)\quad The fuzzy algebra $([0,1],\vee,T)$, where $T$
is a triangular norm,

\noindent (4)\quad $([0,1],\wedge,S)$, where $S$ is a triangular
conorm,

\noindent (5)\quad The tropical semiring
$(\mathbb{R}^{+}\cup\{\infty\},\wedge,+)$, where $\mathbb{R}^{+}$
is the set of all nonnegative real numbers,

\noindent (6)\quad The semiring of all ideals of a commutative
ring with identity together with addition and multiplication of
ideals.
\medskip

\noindent {\bf Lemma 5.1.}\quad  For an incline $R$, the following
hold.
\medskip

\noindent (1)\quad  If $R$ has a maximal element, then it is the
greatest element of $R$.

\noindent (2)\quad  If $R$ has a minimal element, then it is the
least element of $R$.
\medskip

\noindent {\bf Proof.}\quad  (1)\quad  Let $a\in R$ be maximal in
$R$. For any $r\in R$, $a\leqslant r+a$ and so $r+a=a$, i.e.
$r\leqslant a$.

\noindent (2)\quad  Let $a\in R$ be minimal in $R$. For any $r\in
R$, $ra\leqslant a$ and so $a=ra\leqslant r$.\quad $\square$
\medskip

An element in an incline $R$ is the least element of $R$ iff it is
the zero of $R$.
\medskip

\noindent {\bf Lemma 5.2.}\quad  Any incline $R$ with
$|R|\geqslant 3$ is not $k$-simple.
\medskip

\noindent {\bf Proof.}\quad   Since $|R|\geqslant 3$, it follows
from Lemma 5.1 that there exists an element $r\in R$ which is
neither maximal nor minimal in $R$. Let $A=\{x\in R\mid x\leqslant
r\}$. Then $A\neq R$, and $A\neq \{0\}$ if $0\in R$. Show that $A$
is a $k$-ideal of $R$. If $x,y\in A$, then $x\leqslant r$,
$y\leqslant r$ and $x+y\leqslant r$, thus $x+y\in A$. If $x\in A$
and $s\in R$, then $sx\leqslant x\leqslant r$, i.e. $sx\in A$ and
similarly $xs\in A$. If $x\in A$ and $s\in R$ with $s+x\in A$,
then $s\leqslant s+x\leqslant r$ and so $s\in A$. Hence $R$ is not
$k$-simple.\quad $\square$
\medskip

\noindent {\bf Example 5.3.}\quad  There exist exactly two
$2$-element inclines (up to isomorphism). They are $R_0$ and $R_1$
determined as follows.

\[
R_0:
\begin{array}{cc}
\begin{tabular}{c|cc}
+& 0& 1\\
\hline
0& 0& 1\\
1& 1& 1\\
\end{tabular}
&
\begin{tabular}{c|cc}
$\cdot$& 0& 1\\
\hline
0& 0& 0\\
1& 0& 0\\
\end{tabular}
\end{array}
\qquad R_1:
\begin{array}{cc}
\begin{tabular}{c|cc}
+& 0& 1\\
\hline
0& 0& 1\\
1& 1& 1\\
\end{tabular}
&
\begin{tabular}{c|cc}
$\cdot$& 0& 1\\
\hline
0& 0& 0\\
1& 0& 1\\
\end{tabular}
\end{array}
\]
Note that $R_0$ and $R_1$ are congruence-simple and ideal-free.
\medskip

Summarizing Lemmas 5.1 and 5.2 and Example 5.3 gives the following
conclusion.
\medskip

\noindent {\bf Theorem 5.4.}\quad  The only $k$-simple inclines
are $R_0$ and $R_1$.
\medskip

$ $

{\section*{\noindent \bf 6. A counterexample}}

$ $

In this section, we give a counterexample to show that the sum of
$k$-ideals is not necessarily a $k$-ideal in semirings.

Lemma 2.12(i) in Atani \cite{Atani2007} states that if $I$ and $J$
are $k$-ideals of $R$, then $I+J$ is a $k$-ideal of $R$, where $R$
is any commutative semiring. Furthermore, the proofs of Lemma
2.13(i) and Theorem 2.16 use Lemma 2.12(i) in \cite{Atani2007}.

The following example says that Lemma 2.12(i) in \cite{Atani2007}
is false (see also p.75 in \cite{Golan1999} and p.86 in
\cite{Hebisch1998}).
\medskip

\noindent {\bf Example 6.1.}\quad  Let $\mathbb{Z}^+$ be the
semiring of all nonnegative integers together with the usual
addition and multiplication. Though $(2)=2\mathbb{Z}^+$ and
$(3)=3\mathbb{Z}^+$ are $k$-ideals of $\mathbb{Z}^+$, the sum
$(2)+(3)$ is not $k$-closed in $\mathbb{Z}^+$. In fact, both $6$
and $7$ are in $(2)+(3)$ and $1+6=7$, but $1$ is not in $(2)+(3)$.
\medskip

$ $

{\section*{\noindent \bf 7. An open problem}}

$ $

The classification of $k$-simple semirings remains open.
\medskip

\vspace{0.8cm}

{\renewcommand\baselinestretch{1}
\renewcommand\refname}

\end{document}